\documentclass[final]{siamltex}
\usepackage{epsfig}
\usepackage{amsfonts}
\def \abs#1{{\left|#1\right|}}
\def \Norm#1{{\|#1\|}}
\def \inner#1#2{\langle #1, #2 \rangle}
\def \A{\mathcal{A}}   
\def \K{\mathcal{K}}   
\def \V{\mathcal{V}}
   
\def \R{\mathbb{R}}   
\def \Z{\mathbb{Z}}
\def \Sphere{\mathbb{S}}
\def \half{{1\over 2}}
\def \Trace{\mbox{tr}}
\def \E#1{$\times 10^{#1}$}
\renewcommand{\hat}{\widehat}
\renewcommand{\phi}{\varphi}

\title{A spectral method for integral formulations of 
medium-frequency scattering problems} 

\author{Johannes Tausch\footnotemark[2]}

\begin{document}
\maketitle
\footnotetext[2]{Department of Mathematics, Southern Methodist University,
  Dallas, TX 75275, \texttt{tausch@smu.edu}}

\begin{abstract}
  A fast method for the computation of layer potentials that arise in
  acoustic scattering is introduced.  The principal idea is to split
  the singular kernel into a smooth and a local part.  The potential
  due to the smooth part is computed efficiently using non-equispaced
  FFTs, the potential due to the local part is expanded as a series in
  the mollification parameter. The complexity of the approach is shown
  to be $O(n + \kappa^3 \log \kappa)$, where $n$ is the number of
  degrees of freedom in the discretization and $\kappa$ is the
  wave number. The constant factor in this asymptotic estimate is
  small since no singular surface integrals must be computed.
  Therefore the method is particularly efficient for medium-sized
  scatterers (50-100 wavelengths) that may have complicated geometry.
\end{abstract}


\pagestyle{myheadings}
\markboth{J. TAUSCH}{Spectral Integral Method for Scattering Problems}

\section{Introduction}
It is commonly accepted that the boundary element method is an
effective approach to solve the Helmholtz equation in the exterior of
a scatterer. It has the advantage that only the finite boundary
surface has to be discretized and that the radiation condition is
automatically satisfied. In the recent past a variety of methods have
been developed to handle the dense matrices associated with
discretized layer potentials. These can be roughly classified into two
groups, namely hierarchical and grid-based methods.

Examples of hierarchical methods are the Fast Multipole
Method\cite{greengard-rokhlin87}, hierarchical
matrices\cite{borm-grase-hackb02} and wavelets\cite{beyl-coif-rokh91}.
These methods are based on clustering interactions between panels in a
hierarchical manner; the larger the separation, the larger the
clusters. The efficiency and accuracy of these methods depends
critically on how the cluster interactions can be approximated by
low-rank matrices. In the case of boundary integral operators
associated with the Laplace, Stokes or Lam\'e equations asymptotically
optimal schemes have been developed.  That is, the complexity of a
matrix-vector multiplication is order $n$, or order $n \log^p n$,
while the convergence rate of the discretization scheme is preserved,
see, e.g.~\cite{schneider98,sauter00,tausch04}.

In the case of the Helmholtz equation, the size of the scatterer,
measured in wavelengths, is the dominant factor that influences
computational cost and accuracy. It is well known that in the
high-frequency regime large clusters are no longer approximated by
low-rank matrices, and therefore the efficiency of the aforementioned
methods breaks down. To overcome this problem it has been proposed to
use the Fast Multipole Method with diagonal translation operators,
\cite{rokhlin90}. This technique has been extended in 
\cite{darve00a, song-etal98}.

A different approach that avoids large clusters is to replace the
surface distribution by equivalent sources on a uniform grid. The fast
Fourier transform can be employed to compute grid potentials
efficiently. Since the grid is only accurate when the source and the
evaluation point are well separated, the nearby interactions must be
computed directly, by adding up contributions of individual sources.
Grid-based methods are quite popular even though it appears that the
asymptotic complexity is generally higher than what can be achieved
with hierarchical methods. However, in many engineering applications
the geometry is complicated and the mesh is relatively coarse,
therefore constant factors often play an important role.  Applications
of grid-based methods for the Laplace equation can be found in
\cite{hockney-eastwood88, phillips-white97} for elasticity in
\cite{peirce-napier95}. For high-frequency scattering this
methodology, combined with local high-order discretizations have been
described in~\cite{bruno-kunyansky01}. 

If the scatterer is smooth and isomorphic to a sphere, spectral element
methods, based on expansion of the solution into spherical harmonics
have been shown to be successful~\cite{ganesh-graham04}.

The approach described in this paper is closer to grid-based methods
in that FFTs are used to accelerate the matrix-vector product.
However, there is no uniform grid with equivalent charges.  The idea
here is to split the Green's function into a smooth approximation and
a singular, essentially local part. The smooth part of the Green's
function is replaced by a rapidly converging Fourier series.  We will
show how non-equispaced FFTs can be used to compute layer potentials
with such a kernel effectively.  A similar idea has been applied
earlier to the heat equation~\cite{greengard-strain90}.

The local part can be evaluated using expansions with respect to the
mollification parameter. Thus the computation of the local part
amounts to multiplying with a diagonal matrix. Since there is no need
to compute the nearfield directly, we believe that the discussed
approach is competitive with hierarchical and grid-based methods.

We will discuss how the mollification parameter and the number of
Fourier modes have to be selected as a function of the meshwith and
the wave number to obtain efficient and accurate schemes. Our analysis
is based on bounding the error of the bilinear form when the wave
number is increased. It should be noted that this does not give
estimates of the error of the solution. For that, realistic estimates
of constants in the inf-sup condition are necessary, which are not
available. However, we will present numerical examples that suggest 
that the selection strategy of the parameters indeed control the error
when increasing the wavenumber.

\section{Problem Formulation}
For simplicity of exposition, the focus of this
paper will be on the sound-soft acoustic scattering of an incoming field
$u^{\tiny\mbox{inc}}$ of a smooth obstacle $D \subset
\mathbf{R}^3$. The reflected field $u$ is described 
by the Dirichlet problem to the Helmholtz equation with the Sommerfeld
radiation condition
\begin{equation}\label{helm-pde}
\begin{array}{rcl}
\Delta u(x) + \kappa^2 u(x) &=& 0, \quad x \in \R^3 \setminus D\\
u(x) &=& -u^{\tiny\mbox{inc}}(x) , \quad x \in S := \partial D\\
{\partial u \over \partial r} - i \kappa u &=& O({1\over \Norm{x}^2}).
\end{array}
\end{equation}
Here, $\kappa$ is the wave number. We assume
that the problem is scaled such that the scatterer is located inside a
cube of side length $1-d$, that is
\begin{equation}\label{scaling}
S \subset [0,1-d]^3
\end{equation}
where $0<d\ll 1$ is a constant.

A classical approach to treat the Helmholtz problem (\ref{helm-pde})
is the combined layer ansatz of Brakhage and Werner
\cite{brakhage-werner65}, where the scattered field is represented by
a combination of a single- and double layer potential
\begin{equation}\label{combo-layer}
u(x) = (\K - i\eta \V)\sigma(x),\quad x \in \R^3 \setminus \overline{D},
\end{equation} 
where $\eta>0$ is the coupling parameter, $\sigma$ an unknown surface density and 
\begin{eqnarray}
\label{def-single} \V\sigma(x) &:=& \int_S 
    {\exp(i\kappa \Norm{x-y})\over 4\pi\Norm{x-y}}\, \sigma(y)\, dS_y,\\
\label{def-double} \K\sigma(x) &:=& \int_S {\partial \over \partial n_y} 
    {\exp(i\kappa \Norm{x-y})\over 4\pi\Norm{x-y}}\, \sigma(y)\, dS_y
\end{eqnarray}
are the single and double layer operator, respectively. By letting
$x\to S$ from the exterior of the scatterer, and taking the jump relations
of layer potentials under consideration, the following boundary
integral equation for $\sigma$ can be derived
\begin{equation}\label{combined-bem}
{1\over 2}\sigma(x) + (\K - i\eta \V)\sigma(x) = -u^{\tiny\mbox{inc}}(x),
\qquad x \in S.
\end{equation} 
It is well known that (\ref{combined-bem}) is a well posed problem
when $\eta > 0$, see, e.g., \cite{colton-kress83}.

To obtain a discretization of (\ref{combined-bem}) we introduce the
space $X_h$ of piecewise polynomial functions on a triangulation of
$S$. This triangulation is assumed to be quasi-uniform and the maximal
diameter of a triangle is denoted by $h$. The nodal basis $\left\{
  \phi_i \right\}$, $i=1,\dots,n$ of $X_h$ consists of functions with
local support. Note that $n=O(h^{-2})$. The Galerkin approximation is
defined to be the function $\sigma_h \in X_h$ whose residual is
$L^2(S)$-orthogonal to $X_h$. This leads to the linear system $(M+A)x=b$,
where the coefficients of the system matrix and the right hand side
are given, respectively, by
\begin{eqnarray*}
M_{i,j} &=& \half \int_S \phi_i(x) \phi_j(x) dS_x\\
A_{i,j} &=& \int_S\!\!\int_S 
    \left( {\partial \over \partial n_y} - i\eta \right) 
    {\exp(i\kappa \Norm{x-y})\over 4\pi\Norm{x-y}}\, 
    \phi_i(x) \phi_j(y)\, dS_y dS_x\,,\\
b_i &=& -\int_S  \phi_i(x) u^{\tiny\mbox{inc}}(x) dS_x\,.
\end{eqnarray*}
Since this system is large, iterative methods for its solution must be
employed. The dominant cost in such a scheme is the multiplication of
a vector with the dense matrix $A$, which, if it is done in the
obvious way, has $O(n^2)$ complexity. The article will discuss a
scheme to compute the product approximately with a highly reduced flop
count.

\section{Splitting of the Helmholtz Kernel} 
The heart of the method is the splitting of the Helmholtz kernel 
\begin{equation}\label{def-split-kern}
G(r) = G_\delta(r) + E_\delta(r)
\end{equation}
into a smooth part $G_\delta(r)$ and a singular, local part
$E_\delta(r)$. Here, $\delta$ is the mollification parameter that
controls the smoothness of $G_\delta$. This splitting results in a
splitting of the single layer potential
$$
\V g(x) = \Phi^S(x) + \Phi^L(x)
$$
where
\begin{eqnarray}
\Phi^S(x) &=& \int_S G_\delta(x-y) g(y)\, dS_y \label{def-single-sm}\\
\Phi^L(x) &=& \int_S E_\delta(x-y) g(y)\, dS_y \label{def-single-loc}
\end{eqnarray}
The splitting of the double layer potential 
$\K g(x) = \Psi^S(x) + \Psi^L(x)$ is defined analogously.

\subsection{Smooth Part}
The Green's function can be expressed in Fourier space, 
\begin{equation}\label{green-fourier0}
G(r) = {\exp(i \kappa \Norm{r})\over 4\pi\Norm{r}} 
     = {1\over (2\pi)^3 } \int_{\R^3} 
   {1\over \Norm{\omega}^2-\kappa^2} \exp(i r\cdot \omega) \, d^3\omega.
\end{equation}
Exploiting the spherical symmetry of the Fourier transform leads to
\begin{equation}\label{green-fourier}
G(r) = {1\over 2\pi^2}
\int_0^\infty {\rho^2 \over \rho^2-\kappa^2}\,j_0(\rho\Norm{r})\, d\rho
\end{equation}
where $j_0(z) = \sin(z)/z$ is the spherical Bessel function of order
zero.  The integral in (\ref{green-fourier}) is understood in the
sense that the singularity at $\rho = \kappa$ is circumvented in lower
complex half-plane, thereby enforcing the Sommerfeld radiation condition,
see, e.g., \cite{duffy01}.

The decay rate of the transform at infinity determines the regularity
of the kernel. The integrand in (\ref{green-fourier}) is only
$O(\rho^{-1})$ as $\rho \to \infty$ which explains the singularity
of the Green's function in the origin.  A smooth
approximation of the kernel can be obtained by multiplying the
transform with a filter to increase the decay rate at infinity
\begin{equation}\label{green-fourier-sm}
G_\delta(r) = {1\over 2\pi^2} \int_0^\infty 
 H\Big(\delta (\rho^2-\kappa^2) \Big) {\rho^2 \over \rho^2-\kappa^2}
  j_0(\rho\Norm{r}) \, d\rho
\end{equation}
where $H$ is the filter.  Because of the singularity of the integrand
at $\rho=\kappa$ it is more convenient to write the filter in the form
as it appears in (\ref{green-fourier-sm}) and not as $H(\delta
\rho^2)$.

There are several possible choices for $H$. If the filter is a
rational function, then $G_\delta$ can be expressed in closed form. To
that end, write the filter in partial fraction decomposition
$$
H(z) = \sum_{k=1}^q {c_k\over z + w_k^2}\,,
$$
where the coefficients $w_k$ and $c_k$ are at our disposition. 
It will become clear later that 
because of the singularity of the integrand it is necessary that
\begin{equation}\label{chi-one}
H(0) = 1.
\end{equation}
Basic complex variable arguments show that
\begin{equation}\label{filter-pfd}
{H(z)\over z} = \sum_{k=0}^q {d_k\over z + w_k^2}\,.
\end{equation}
where $w_0 = 0$ and because of condition (\ref{chi-one})
\begin{equation}\label{def-dk}
d_0 = 1 \quad\mbox{and}\quad  d_k = {c_k\over w_k^2}\,.
\end{equation}
Substitution of (\ref{filter-pfd}) into (\ref{green-fourier-sm}) and
the change variables $\rho \to \sqrt{\delta} \rho$ leads to
$$
G_\delta(r) = {1 \over 2\pi^2 \sqrt{\delta}}
  \sum_{k=0}^q d_k\, \int_0^\infty \!
  { \rho^2 \over \rho^2-\hat w_k^2}\, 
   j_0\left(\rho {\Norm{r}\over\sqrt{\delta}}\right)\, d\rho,
$$
where 
\begin{equation}\label{def-wk}
\tilde \kappa = \sqrt{\delta} \kappa  \quad\mbox{and}\quad
\tilde w_k = \sqrt{ \tilde \kappa^2 - w_k^2}
\end{equation}
The integrals in the last expression are of the same form as
(\ref{green-fourier}). Therefore the mollified Green's function has
the closed form
\begin{equation}\label{decomp-1}
G_\delta(r) = {\exp( i \kappa \Norm{r})\over 4\pi \Norm{r}} +
  \sum_{k=1}^q d_k {\exp( i \tilde w_k \Norm{r}/\sqrt{\delta})\over 4\pi \Norm{r}}
\end{equation}
In the discussion below, it will be convenient to write the
decomposition in (\ref{decomp-1}) in the form
\begin{equation}\label{decomp-2}
G(r) = G_\delta(r) +
   {1\over\sqrt{\delta}} E\left(\Norm{r}\over \sqrt{\delta}\right)
\end{equation}
where $E$ is singular at $z=0$ given by
\begin{equation}\label{def-remain}
E(z) = - \sum_{k=1}^q d_k {\exp( i \tilde w_k z)\over 4\pi z}\,.
\end{equation}
A good filter must satisfy two properties. First, it must decay
rapidly to ensure smoothness, that is, there must be
a constant $c$ such that
\begin{equation}\label{decay-filter}
\abs{H(z)} \leq c \min\left( 1, z^{-q} \right) \,.
\end{equation}
Second, the kernel $E$ must decay exponentially away from the origin,
this is why it will be referred to as the local part. To ensure the latter
condition, it is necessary that 
Im$(\tilde w_k)$ is bounded away from zero as $\kappa \to \infty$.
If $w_k$ is real, then (\ref{def-wk}) implies that
\begin{equation}\label{delta-kappa}
\sqrt{\delta} < {1\over \kappa} \min_{k>0} \abs{w_k}.
\end{equation}
The latter condition implies that $\tilde \kappa$ is bounded as
$\kappa \to \infty$.

An example of a filter that satisfies
(\ref{decay-filter}) is given by
\begin{equation}\label{filter1}
H(z) = \prod_{k=1}^{q} {1\over k+z},
\end{equation}
This filter has poles $w_k = i\sqrt k,\;k=1,\dots,q$. Condition 
(\ref{delta-kappa}) is equivalent to $\sqrt{\delta} < 1/\kappa$.
The filter
\begin{equation}\label{filter2}
H(z) = {1\over (1+z)^q}
\end{equation}
also satisfies the decay property. Because of the repeated
poles, the smooth part corresponding to this filter is slightly
different from (\ref{decomp-1}): 
$$
G_\delta(r) = {\exp( i \kappa \Norm{r})\over 4\pi \Norm{r}} +
  p\left({\Norm{r}\over\sqrt\delta}\right) 
  {\exp( i \tilde w \Norm{r}/\sqrt\delta) \over \Norm{r}}.
$$
Here $p(\cdot)$ is a polynomial of degree $q-1$ and $\tilde w =
\sqrt{1-\delta \kappa^2}$. As with the previous filter, condition
(\ref{delta-kappa}) is $\sqrt{\delta} < 1/\kappa$.

\subsection{Local Part}
The smooth part is
a good approximation of the actual Green's function if $\delta$ is
small and $r$ is large. In the neighborhood of the origin the two
functions are very different and therefore the contribution of the
local part must be accounted for. In this section
we show that the local part has an expansion with respect to the
mollification parameter $\sqrt{\delta}$ and show how to compute the
expansion coefficients. 

Since condition (\ref{delta-kappa}) implies that the coefficients
$\tilde w_k$ in (\ref{def-remain}) have a positive imaginary part, the
function $E_\delta$ decays exponentially away from the origin. We
introduce the smooth cut-off function $\tilde \chi$ for some $0<\nu<1$ which
is small enough such that the surface has a parameterization of the
form $y(t) = x + At + nh(t)$ in the $\nu$-neighborhood of $x$. Here
$n$ is the normal of the surface at the point $x$, $A \in \R^{3\times 2}$ 
has two orthogonal columns that span the tangent plane at $x$
and $h(t) = O(|t|^2)$ is some scalar function in $t \in \R^2$.  The
local single-layer potential $\Phi_\delta(x)$ in
(\ref{def-single-loc}) can be written in the form
\begin{eqnarray}
  \Phi_\delta(x) &=& \int_{S} E_\delta(x-y) g(y)\, dS_y \nonumber \\
    &=& \int_{S} E_\delta(x-y) \tilde \chi_\nu(x-y) g(y)\, dS_y + 
       O\left(\exp\left(-{\nu\over\delta}\right) \right) \nonumber \\
    &=& \int_{\R^2} E_\delta(t) \tilde g(t) d^2t
    + O\left(\exp\left(-{\nu\over\delta}\right) \right).
\label{expan-phi}
\end{eqnarray}
Here, $E_\delta(t) = E_\delta(x-y(t))$, 
$\tilde g(t) = \tilde \chi_\nu(x-y(t)) g(t) J(t)$ and $J(t)$ is the Jacobian of
the parameterization. For simplicity of the argument we assume that
the function $h(t)$ in the parameterization of the surface is
analytic, that is, 
\begin{equation}\label{expan-h}
h(t) = \sum_{|\alpha| \geq 2} h_\alpha t^\alpha.
\end{equation}
Thus there are there are $C^{\infty}$-functions $H_n$ such that
\begin{equation}\label{expan-r}
r(t) := \Norm{x-y(t)} = \Norm{t} \sum_{n=0}^\infty \Norm{t}^n H_n(\hat t) 
\end{equation}
where $\hat t := t/|t|$ and $H_0(\hat t)=1$ and $H_1(\hat t)=0$. 
Substituting (\ref{expan-r}) into (\ref{expan-phi}) leads to
\begin{eqnarray}
  \Phi_\delta(x) &=& {1\over \sqrt{\delta}} \int_{\R^2} 
  E\left({r(t)\over \sqrt{\delta}}\right) \tilde g(t) \, d^2t +
       O\left(\exp\left(-{\nu\over\delta}\right) \right) \\
  &=& \sqrt{\delta} \int_{\R^2} 
  E\left(|t| \sum_{n=0}^\infty (\sqrt{\delta} |t|)^n H_n(\hat t) \right) 
  \tilde g(\sqrt{\delta} t) \, d^2t +
  O\left(\exp\left(-{\nu\over\delta}\right) \right) 
\end{eqnarray}
where the second integral is the result of the change of variables
$t\mapsto t/\sqrt{\delta}$. The integral as a function of
$\sqrt{\delta}$ is $C^\infty$, and can therefore be expanded in a
Taylor series. The expansion coefficients are derivatives of the integral
with respect to $\sqrt{\delta}$. We see that
\begin{equation}\label{expan-phi2}
\Phi_\delta(x) = \delta^{1\over 2} \Phi_0 g(x) + O(\delta^{3\over 2})
\end{equation}
where
$$
\Phi_0 = \int_{\R^2} E(\Norm{t}) d^2t = 
{\pi\over i} \sum_{k=1}^q {d_k\over w_k}\,.
$$
The double layer potential is given by
$$
\Psi_\delta(x) = -{1\over\delta} 
   \int_S E'\left({\Norm{x-y}\over \sqrt\delta}\right) 
   {(x-y)\cdot n_y \over \Norm{x-y} } g(y) \,dS_y .
$$
The second factor of the kernel can be expanded in a similar manner as
(\ref{expan-r}), we find that
$$
{(x-y)\cdot n_y \over \Norm{x-y} } =
-\Norm{t} \left( h_{02} \cos^2\theta
+ h_{11} \cos\theta\sin\theta + h_{20} \sin^2\theta \right) 
+ O(\Norm{t}^2)
$$
where $\theta$ is the angular coordinate of $t$ and the $h$'s are from
expansion (\ref{expan-h}). Proceeding in a similar manner as for the
single layer operator the following expansion follows
\begin{equation}\label{expan-psi2}
\Psi_\delta(x) = \delta^{1\over 2} \Psi_0 g(x) + O(\delta^{3\over 2})
\end{equation}
where
$$
\Psi_0 = -{\pi\over i} (h_{02} + h_{20}) \sum_{k=1}^q {d_k\over w_k}\,.
$$

\section{A fast algorithm for smooth, periodic kernels}
We describe a fast algorithm for the smooth part of the single layer in
(\ref{def-single-sm}), which is based on Fourier
analysis. Modifications for the double layer are minimal and mentioned
at the end of the section. For rapid convergence of the Fourier series
it is necessary to multiply the smooth kernel with a 
sufficiently smooth cut-off function $G^S := \chi G_\delta$ 
that is unity inside the cube $[-1+d,1-d]^3$ and vanishes outside
$[-1,1]^3$. Recall that we assumed in (\ref{scaling}) that the surface is
contained in $[0,1-d]^3$, thus the cut-off function has no effect in
the integral, and the smooth part is given by
\begin{equation}\label{def-phi}
\Phi^S(x) = \int_S G^S(x-y) g(y) \, dS_y, \qquad x\in S.
\end{equation}
The kernel $G^S$ can be approximated by the truncated Fourier
series $G_N$
\begin{equation}\label{fourier-series-gf}
G_N(r) := \sum_{\Norm{k}_\infty \leq N} \hat G_k \exp(\pi i\, k^T r)\,,
\quad r \in [-1,1]^3
\end{equation}
where the summation index $k$ is in $\Z^3$. 
The resulting approximate potential is given by
\begin{equation}\label{fourier-series-phi}
\Phi_N(x) = \int_S G_N (x-y) g(y) dS_y 
   = \sum_{\Norm{k}_\infty \leq N} \exp(\pi i k\cdot x) \hat d_k
\end{equation}
where $\hat d_k = \hat G_k \hat g_k$ and 
\begin{eqnarray}
\hat G_k &=& {1\over 8} \int_{[-1,1]^3} \exp(-\pi i k\cdot r) G^S(r)\, d^3r,
\label{fourier-coeff-G}\\
\hat g_k &=& \int_S \exp(-\pi i k\cdot y) g(y) dS_y.
\label{fourier-coeff-g}
\end{eqnarray}
In case the smooth part of the double layer is to be calculated, the
coefficients $\hat g_k$ must be replaced by
\begin{equation}\label{fourier-coeff-g2}
\hat g_k = \int_S {\partial\over \partial n_y} \exp(-\pi i k\cdot y) g(y) dS_y.
\end{equation}
To simplify the discussion, our notations will not distinguish between
the coefficients in (\ref{fourier-coeff-g}) and (\ref{fourier-coeff-g2}).
In summary, the potential computation of the potential due to the
smooth parts consists of three stages.
in~(\ref{fourier-series-phi}).
\begin{enumerate}
  \item Compute the Fourier coefficients $\hat g_k$ in~(\ref{fourier-coeff-g}).
  \item Multiply $\hat d_k := \hat G_k \hat g_k$ for $\Norm{k}_\infty\leq N$.
  \item Evaluate the Fourier series~(\ref{fourier-series-phi}) for $x\in S$.
\end{enumerate}
The choice of the truncation parameter $N$ depends on the wave number and the
mollification parameter and can be much smaller than the linear system size
$n$. In Section~\ref{sec-error} the exact dependence will be
investigated.  Stage 2 obviously involves $O(N^3)$ operations, the
other two stages can be executed efficiently using non equispaced Fast
Fourier Transforms. This will be discussed next.

\subsection{Computation of the $\hat g_k$'s}\label{sec-comp-gk}
\begin{figure}[hbt]
\begin{center}
\epsfig{figure=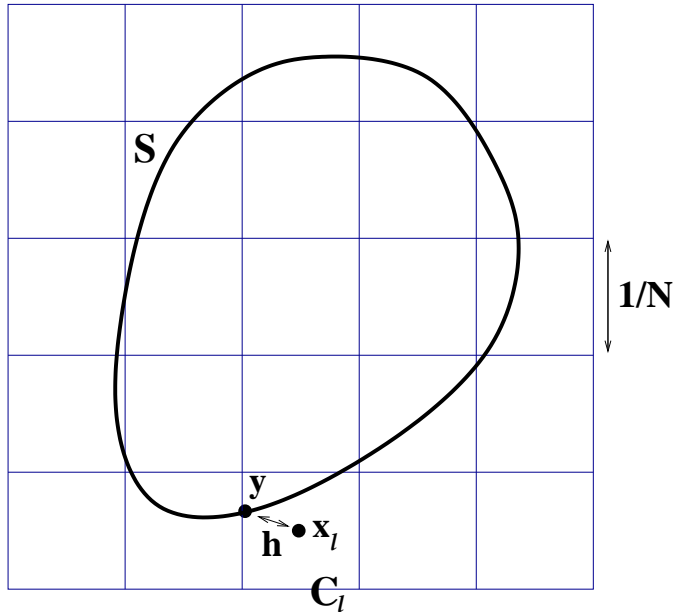,width=3.5in}
\end{center}
\caption{Two dimensional illustration of the geometry.}
\label{fig-grid}
\end{figure}

In this section we describe how FFTs can be used to efficiently
compute the Fourier coefficients of the function $g$. To that end, the
three-space is divided into small cubes $C_l$, $l=(l_1,l_2,l_3)\in
Z^3, 0\leq l_j <N$. These cubes have centers $x_l = l/N$ and
side length $1/N$. Note that $N$ is the same as in
(\ref{fourier-series-gf}) and therefore the cubes get smaller if more
terms in the Fourier series expansion of the Green's function are
retained.  Because of assumption (\ref{scaling}) $S$ is contained in
the union of all cubes and set $S_l = C_l \cap S$ to denote the piece
of the surface that intersects with the $l$th cube, c.f.\ 
Figure~\ref{fig-grid}.

From (\ref{fourier-coeff-g}) it follows that
the Fourier coefficients of $g$ can be written as
\begin{equation}\label{fourier-sum-g}
\hat g_k = \sum_{\Norm{l}_\infty \leq N} 
     \exp\left(-\pi i k\cdot l\over N\right) 
     \int_{S_l} \exp(-\pi i k\cdot(y-x_l) ) g(y)\,dS_y.
\end{equation}
The frequency and the spatial variable in the integral can be
separated using the Jacobi-Anger expansion
$$
\exp( -i \xi t ) = \sum_{\nu=0}^\infty (-i)^\nu (2\nu + 1)
j_\nu(\xi) P_\nu(t), \qquad -1 \leq t \leq 1, 
$$
see, e.g., \cite{nedelec01}.  Here, $j_\nu(\cdot)$ is the spherical
Bessel function of order $\nu$ and $P_\nu(\cdot)$ is the Legendre
polynomial of degree $\nu$. This formula generalizes to the
three-variate case and can be applied to the
integrand in (\ref{fourier-sum-g})
\begin{equation}\label{trunc-ja}
\exp( -\pi i k\cdot(y-x_l) ) \approx 
  \sum_{|\alpha| \leq p} (-i)^{|\alpha|} (2\alpha+1) 
      j_\alpha(\pi k H) P_\alpha({y-x_l\over H})
\end{equation}
where $p$ is the expansion order, $H = 1/(2N)$,
$\alpha=(\alpha_1,\alpha_2,\alpha_3)$ is a multi-index, 
$|\alpha|=\alpha_1+\alpha_2+\alpha_3$, 
$j_\alpha(x)=j_{\alpha_1}(x_1)j_{\alpha_2}(x_2)j_{\alpha_3}(x_3)$ and
$P_\alpha(x)$ is defined similarly. Substitution of (\ref{trunc-ja}) 
into (\ref{fourier-sum-g}) leads to the approximation
$$
\hat g_k \approx \sum_{|\alpha|\leq p} (-i)^\alpha (2\alpha+1)
j_\alpha(\pi k H) \sum_{\Norm{l}_\infty \leq N} 
 \exp\left(-\pi i k\cdot l\over N\right) m_l^\alpha(g),
$$
where
\begin{equation}\label{mom-surf}
  m_l^\alpha(g) = \int_{S_l} P_\alpha({y-x_l\over H}) g(y)\,dS_y
\end{equation}
is a moment for which exact formulas can be derived if the function
and the surface are discretized. In particular, if $g$ is a piecewise
polynomial, then the moments are linear transformations of the
coefficients of $g$ corresponding to the nodal basis. The matrix that
maps the coefficients to the $\alpha$-th moments is denoted by
$M_\alpha$. The number of nonzero entries in $M_\alpha$ is $n$. 

In matrix form, the (approximate) coefficient vector $\hat g$ is given by
\begin{equation}\label{g-fft}
\hat g = \sum_{|\alpha|\leq p} K_\alpha F\, M_\alpha \vec{g},
\end{equation}
where $F$ is the $2N$-long three-dimensional discrete Fourier
transform, $\vec{g}$ the vector of coefficients of $g$ and $K_\alpha$
is a diagonal matrix with the factors $(-i)^{|\alpha|} (2\alpha+1)
j_\alpha(\pi k H)$. The computation of $\hat g$ involves
$(p+1)(p+2)(p+3)/6$ FFTs.  In Section~\ref{neFFT} we will show that
it suffices to use a small value of $p$.

\subsection{Evaluation the Fourier series} 
In the Galerkin discretization, the $i$-th component of the matrix-vector
product $\Phi_i$ is the inner product of the potential $\Phi$
in~(\ref{def-phi}) with the $i$-th nodal basis function $\varphi_i$.
For the fast method, the potential is replaced with the approximated
potential $\Phi_N$ in~(\ref{fourier-series-phi}). In order to evaluate
the potential efficiently, the Jacobi-Anger approximation
(\ref{trunc-ja}) is used again, in a very similar manner as in the
previous section. This is shown in the following computation
\begin{eqnarray*}
\Phi_i &=& \int_S \varphi_i(x) \Phi_N(x) \, dS_x\\
   &=& \sum_{\Norm{k}_\infty\leq N} \exp\left(\pi i\, k\cdot l\over N\right) 
     \int_S \exp(\pi i k\cdot(y-x_l) ) \varphi_i(x) \, dS_x\, d_k\\
   &\approx& \sum_{|\alpha|\leq p} \sum_{\Norm{k}_\infty \leq N}
         \exp\left(\pi i\, k\cdot l\over N\right) 
         i^{|\alpha|} (2\alpha+1) j_\alpha(\pi k H)
         m_l^\alpha(\varphi_i) \, d_k \,.
\end{eqnarray*}
In matrix notation, the above can be written as
\begin{equation}\label{phi-fft}
\vec \Phi = \sum_{|\alpha|\leq p} M_\alpha^T F^* K_\alpha \vec d.
\end{equation}
Hence $(p+1)(p+2)(p+3)/6$ FFTs are necessary to compute the vector
$\vec \Phi$. Furthermore, it is evident that the operation
(\ref{phi-fft}) is the adjoint of operation (\ref{g-fft}).

\section{Error Analysis}\label{sec-error}
In this section we derive estimates for the error of the bilinear form
introduced when the mollified kernel corresponding to the single layer
operator is replaced by the truncated Fourier series expansion. To
keep the technical level of the discussion at a minimum we omit the
discussion of the double layer operator, because it is completely
analogous to the single layer. The main concern is the situation where
$\kappa\to \infty$ and our goal is to determine $N$ and $p$ as a
function of the wavenumber such that the resulting error remains
bounded.

It is straightforward to see that
$$
\inner{f}{(\A - \A_N) g} = 
  \sum_{\Norm{k}_\infty>N} \hat G_k \hat f_k \hat g_k
$$
where $\A$ is the surface integral operator with kernel $G^S$, $\A_N$
its Fourier series approximation, $\hat G_k$ are the Fourier
coefficients of $G$ and $\hat f_k$, $\hat g_k$ are the Fourier
coefficients of surface distributions as defined in (\ref{fourier-coeff-g}).
The obvious way to estimate the error is
\begin{eqnarray}
\abs{\inner{f}{(\A - \A_N) g}} &\leq& 
  \sum_{\Norm{k}_\infty > N} \abs{\hat G_k}  
   \sup_k \abs{\hat f_k}
   \sup_k \abs{\hat g_k} \nonumber\\
&\leq& 
 \sum_{\Norm{k}_\infty > N} \abs{\hat G_k} \Norm{f}_{L^2(S)} \Norm{g}_{L^2(S)}.
\label{pre-est-inner}
\end{eqnarray}
It is not possible to work with $l_2$-estimates of $\hat f_k$, because
$f$ can be regarded as surface-delta function in $\R^3$ which cannot
be bounded by the $L^2$-norm. Thus the Fourier coefficients of the
kernel must be estimated in the $l_1$-norm which amounts to an estimate
in the $L^\infty_{[-1,1]^3}$-norm.  

The Fourier coefficients can be related to the derivatives of the
function with the standard integration-by-parts argument. Since the
kernel is a spherically symmetric, three-variate function, it is
convenient to work with the Laplacian and the Green's formula. Because
of $\Delta \exp(ik\cdot r) = -\Norm{k}^2 \exp(ik\cdot r)$ it follows
that
\begin{eqnarray*}
\lefteqn{\hat G_k = 
    {1\over 8} \int_{[-1,1]^3} G_\delta(r) \chi(r) \, d^3 r}\\&&
  = {1\over 8\Norm{k}^{2m}} \int_{[-1,1]^3} \Delta^m(G_\delta(r) \chi(r)) \, d^3 r
  = {1\over \Norm{k}^{2m}} \hat{[\Delta^m(G_\delta \chi)]}_k\,,
\end{eqnarray*}
for any integer $m$ for which the right-hand side is defined. With
this estimate at hand, one obtains
\begin{eqnarray}\label{pre-est-err}
\sum_{\Norm{k}_\infty > N}  \abs{\hat G_k} 
&\leq& \left( 
\sum_{k\in \Z^3} \abs{\hat{[\Delta^m(G_\delta \chi)]}_k}^2 \right)^\half
\left( \sum_{\Norm{k}_\infty > N} {1\over \Norm{k}^{4m}}\right)^\half
\nonumber\\
&\leq&
c N^{{3\over 2}-2m} \Norm{\Delta^m (G_\delta \chi)}_{L^2} 
\end{eqnarray}
where the first step follows from the Cauchy-Schwarz inequality and
the second step follows from Parseval's equation. Since the mollified
Green's function gets more peaked in the origin as $\delta\to 0$ the
norm of $G_\delta$ cannot be treated as a constant. Therefore $N$ must be linked
to $\delta$.  Unfortunately, the right-hand side in
(\ref{pre-est-err}) involves the product of the Green's function with
the cut-off function. The product rule leads to
estimates that involve factors which depend on $m$ and which are difficult
to control. Furthermore, the argument in (\ref{pre-est-err}) assumes
that the cut-off function has the same regularity as $G_\delta$,
which, as numerical experiments suggest, is not necessary.

The following discussion presents a refined analysis intended to obtain more
realistic error estimates.

\subsection{Estimates of the Derivatives}
Derivatives of the mollified Green's function can be obtained either
from the Fourier transform (\ref{green-fourier-sm}) or the closed form
(\ref{decomp-1}). The Fourier
integral leads to estimates that display the dependence on the order
of differentiation more clearly. For the
subsequent error analysis it suffices to work with powers of the
Laplacian. Since the Bessel function satisfies the Helmholtz equation
it follows from (\ref{green-fourier-sm}) that
\begin{equation}\label{green-fourier-D}
\Delta^m G_\delta(r) = (-1)^m \int_0^\infty 
 H\Big(\delta (\rho^2-\kappa^2) \Big) {\rho^{2+m} \over \rho^2-\kappa^2}
  j_0(\rho\Norm{r}) \, d\rho
\end{equation}
for any integer $m < q$. Because of the assumption (\ref{decay-filter}), the filter
decays only if the argument is larger than unity.
This condition leads to 
\begin{equation}\label{def-rhonull}
\delta ( \rho^2 - \kappa^2 ) \geq 1 \Rightarrow 
\rho \geq \left( {1\over \delta} + \kappa^2 \right)^\half =: \rho_0.
\end{equation}
The definition of $\rho_0$ immediately implies that
\begin{equation}\label{def-rhonull2}
\rho_0 = {1\over \sqrt{\delta}} \sqrt{1 + \tilde \kappa^2},
\end{equation}
i.e., $\rho_0 \leq c \kappa$.
\begin{lemma}
The following estimates hold
\begin{eqnarray}
\abs{ \Delta^m G_\delta(r) } \leq 
   c {m \over \Norm{r} } \left( {1+\tilde \kappa^2\over \delta} \right)^m,
   \label{est-G-inf}\\
\Norm{ \Delta^m G_\delta(r) }_{L^2_{[-1,1]^3}} \leq 
   c m \left( {1+\tilde \kappa^2 \over \delta}\right)^m.
   \label{est-G-two}
\end{eqnarray}
\end{lemma}

\begin{proof}
Break integral (\ref{green-fourier-D}) into three parts
$$
\Delta^m G_\delta(r) = {(-1)^m\over 2\pi^2}  \big( I_1 + I_2 + I_3 \big),
$$
where
\begin{eqnarray*}
I_1 &=& \int_0^{\rho_0} {f(\rho)-f(\kappa)\over \rho-\kappa} 
    {\rho^2 \over \rho+\kappa} j_0(\rho\Norm{r}) \, d\rho,\\
I_2 &=& f(\kappa) \int_0^{\rho_0} 
    {\rho^2 \over \rho^2-\kappa^2} j_0(\rho\Norm{r}) \, d\rho,\\
I_3 &=& \int_{\rho_0}^\infty f(\rho) {\rho^2 \over \rho^2-\kappa^2}
   j_0(\rho\Norm{r})\, d\rho,
\end{eqnarray*}
and
$$
f(\rho) := \rho^{2m} H\Big(\delta (\rho^2-\kappa^2) \Big).
$$
Since $\delta \rho_0^2 \leq c$ the derivative of $f(\rho)$ is
bounded by $\abs{f'(\rho)} \leq c m \rho_0^{2m-1}$ for
$0\leq\rho\leq\rho_0$. Thus the first integral can be estimated as follows
\begin{eqnarray}
\abs{I_1} 
&\leq& 
   \max_{0\leq\rho\leq\rho_0} \abs{f'(\rho)} 
   \int_0^{\rho_0}{\rho^2 \over \rho+\kappa}\, \abs{j_0(\rho\Norm{r})}\,d\rho
   \nonumber\\
&\leq& 
   c m {\rho_0^{2m-1}\over r} \int_0^{\rho_0}{\rho \over \rho+\kappa} \,d\rho
   \nonumber\\
&\leq& 
   c m {\rho_0^{2m}\over \Norm{r}}. \label{bound-I1}
\end{eqnarray}
The integral in $I_2$ can be computed in closed form. The resulting
expression involves the integral sine and cosine functions and are
easily shown to be uniformly bounded, thus 
\begin{equation}\label{bound-I2}
\abs{I_2} \leq c {\kappa^{2m}\over \Norm{r}}.
\end{equation}
The integration in $I_3$ is over the interval where the filter is
decreasing. Therefore
\begin{eqnarray*}
\abs{I_3} 
&\leq& 
  {1\over \Norm{r} \delta^q} \int_{\rho_0}^\infty 
      {\rho^{2m+1} \over (\rho^2-\kappa^2)^q } {1\over\rho^2-\kappa^2}\,
      d\rho\\
&\leq& 
  {1\over \Norm{r} \delta^q} 
    \max_{\rho_0 \leq \rho} \abs{\rho^{2m+1} \over (\rho^2-\kappa^2)^q}
     \int_{\rho_0}^\infty {1\over\rho^2-\kappa^2}\, d\rho.
\end{eqnarray*}
The function to be maximized is monotonically decreasing; the integral
can be computed in closed form and estimated by $c/\kappa$. Thus $I_3$
can be estimated by
\begin{equation}\label{bound-I3}
\abs{I_3} \leq {c\over \kappa \Norm{r}} \rho_0^{2m+1}.
\end{equation}
Estimate (\ref{est-G-inf}) is immediate from (\ref{def-rhonull2}),
(\ref{bound-I1}), (\ref{bound-I2}) and (\ref{bound-I3}). Estimate
(\ref{est-G-two}) follows from (\ref{est-G-inf}) because the $1/r$
singularity cancels upon integration.
\end{proof}

Using very similar arguments as in the previous proof, the first and
second derivatives of powers of the Laplacian can be estimated. We
only state the result.
\begin{lemma}
For $\abs{\alpha} \leq 2$ we have
\begin{eqnarray}
\Norm{ \partial^\alpha \Delta^m G_\delta(r) }_{L^2_{[-1,1]^3}} \leq 
   c {m \over \delta^{m+\half} } \left( 1+\tilde \kappa^2
   \right)^{m+{\abs{\alpha}\over 2}} .
   \label{est-gradG-two}
\end{eqnarray}
\end{lemma}

\subsection{Approximation Analysis of the Fourier Series}
Our goal is an estimate the Fourier truncation error in the spirit of 
(\ref{pre-est-err}), that does not involve high-order derivatives of
the cut-off function. 

\begin{lemma}
If $\chi \in C^3$ then the Fourier coefficients of the kernel are given by
\begin{eqnarray*}
\lefteqn{\hat G_k = {1\over \Norm{k}^{2m}}
\left( \hat{[\Delta^m G_\delta]}_k + 
       2m \hat{[\nabla \Delta^{m-1} G_\delta\cdot\nabla\chi]}_k\right)
   + }\\
&&
\sum_{l=0}^{m-1} {1\over \Norm{k}^{2l+2}}
\sum_{\abs{\alpha} \in \{1,2\} \atop \abs{\beta} = 4-\abs{\alpha}}
a^l_{\alpha,\beta} 
\hat{\left[ \partial^\alpha \Delta^{l-1} G_\delta \partial^\beta \chi \right]}_k
\end{eqnarray*}
for any integer $m$ for which the right-hand side is defined. The
coefficients satisfy $a^l_{\alpha,\beta} \in \{0, 1, 2l, 4l\}$.
\end{lemma}

Note that the above expression only contains derivatives of the
cut-off function up to order three. Furthermore, the remainder
(i.e., the sum over $l$), contains derivatives of $G_\delta$ which are
at least two orders lower than the power of $1/\Norm{k}$. This will be
essential for the subsequent error analysis.

\begin{proof}
A simple application of the product rule shows that
\begin{eqnarray*}
\Delta(G_\delta\chi) &=& 
  \Delta G_\delta \chi + 2 \nabla G_\delta\cdot\nabla\chi + G_\delta \Delta \chi,\\
\Delta(\nabla G_\delta\cdot\nabla\chi) &=&
   \nabla \Delta G_\delta\cdot \Delta \chi + 2 \Trace(G_\delta''\chi'') + 
    \nabla \chi \cdot \nabla \Delta G_\delta.
\end{eqnarray*}
From integration by parts it follows that
$$
\hat G_k = {1\over \Norm{k}^2} \hat{[\Delta(G_\delta\chi)]}_k
= {1\over \Norm{k}^2} \left(
\hat{[\Delta G_\delta \chi]}_k
+ 2 \hat{[\nabla G_\delta\cdot\nabla\chi]}_k
+ \hat{[G_\delta \Delta \chi]}_k \right).
$$
Repeating this argument for the first two terms and leaving the third term
unchanged leads to
\begin{eqnarray*}
\hat G_k &=& {1\over\Norm{k}^4} \left(
\hat{[\Delta^2 G_\delta \chi]}_k
+ 4 \hat{[\nabla \Delta G_\delta\cdot\nabla\chi]}_k \right) \\
&+&
{1\over\Norm{k}^4}  \left( 
\hat{[\Delta G_\delta \Delta \chi]}_k
 + 4 \hat{[\Trace G''\chi'']}_k
 + 2 \hat{[\nabla G_\delta \cdot \nabla \Delta \chi]}_k
\right)
+ {1\over\Norm{k}^2}
 \hat{[G_\delta \Delta \chi]}_k.
\end{eqnarray*}
By induction one finds that
\begin{eqnarray*}
\hat G_k &=& {1\over\Norm{k}^{2m}} \left(
\hat{[\Delta^m G_\delta \chi]}_k
+ 2m \hat{[\nabla \Delta^{m-1} G_\delta\cdot\nabla\chi]}_k \right) \\
&+&
\sum_{l=0}^{m-1}
{1\over\Norm{k}^{2l+2}}  \left( 
\hat{[\Delta^l G_\delta \Delta \chi]}_k
 + 4l \hat{[\Trace \Delta^{l-1} G''\chi'']}_k
 + 2l \hat{[\nabla \Delta^{l-1} G_\delta \cdot \nabla \Delta \chi]}_k
\right)
\end{eqnarray*}
which is the assertion.
\end{proof}

Combining this result with the estimates of the derivatives of
$G_\delta$ leads to the next theorem
\begin{theorem}
If $\chi \in C^3$, $f,g \in L^2(S)$ and $N$ is chosen such that 
\begin{equation}\label{def-lambda}
\lambda := \left({1 + \tilde \kappa^2 \over \delta N^2}\right)^\half < 1,
\end{equation}
then for any integer $0<m<q$ the approximation error of the truncated
Fourier series is bounded by
$$
\abs{\inner{f}{(\A-\A_N)g}}
\leq c \left( m^2 N^{3\over 2} \lambda^m  + N^{-\half} \right) 
\Norm{f}_{L^2(S)} \Norm{g}_{L^2(S)}.
$$
\end{theorem}

This theorem suggests how $N$ must be selected to control the error as
$\kappa\to \infty$. Recall that condition (\ref{delta-kappa}) implies
that $\delta\sim 1/\kappa^2$ to ensure that $\tilde \kappa$ is
bounded. Because of (\ref{def-lambda}) one has to select $N$ such that
$N\sim 1/\sqrt{\delta}$. Since $q$ is free, the product $N^{3\over 2}
\lambda^m$ can always be controlled by letting $q$ increase as $N$
increases. The influence of $q$ on the computational cost is
negligible, therefore the error can be controlled with complexity $O(n
+ \kappa^3 \log \kappa)$.

\begin{proof}
Using the previous lemma it follows that
\begin{eqnarray*}
\sum_{\Norm{k}_\infty > N} \abs{\hat G_k} &\leq& 
\sum_{\Norm{k}_\infty > N} {1\over \Norm{k}^{2m}} 
    \abs{\hat{[\Delta^m G_\delta]}_k + 
       2m \hat{[\nabla \Delta^{m-1} G_\delta\cdot\nabla\chi]}_k} \\
&+& 
\sum_{l=0}^{m-1} 
\sum_{\abs{\alpha} \in \{1,2\} \atop \abs{\beta} = 4-\abs{\alpha}}
\sum_{\Norm{k}_\infty > N} {1\over \Norm{k}^{2l+2}}
\abs{a^l_{\alpha,\beta}}
\abs{\hat{\left[ \partial^\alpha \Delta^{l-1} G_\delta \partial^\beta \chi \right]}_k}.
\end{eqnarray*}
Using Cauchy-Schwarz and Parseval in a similar manner that lead to
estimate (\ref{pre-est-err}) 
\begin{eqnarray*}
\sum_{\Norm{k}_\infty > N} \abs{\hat G_k} &\leq& 
{1\over N^{2m-{3\over 2}}} \left( \Norm{\Delta^m G_\delta \chi} +
  \Norm{\nabla \Delta^{m-1} G_\delta\cdot\nabla\chi} \right)\\
&&+\;
\sum_{l=0}^{m-1} {1\over N^{2l+\half}}
\sum_{\abs{\alpha} \in \{1,2\} \atop \abs{\beta} = 4-\abs{\alpha}}
\abs{a^l_{\alpha,\beta}} 
  \Norm{\partial^\alpha \Delta^{l-1} G_\delta \partial^\beta \chi}.
\end{eqnarray*}
The derivatives in the above expression can be estimated using the
inequalities (\ref{est-G-two}) and (\ref{est-gradG-two}) derived in
the previous section. Recalling the definition of $\lambda$ in
(\ref{def-lambda}), this leads to
$$
\sum_{\Norm{k}_\infty > N} \abs{\hat G_k} \leq 
c N^{3\over 2} m^2 \lambda^m +  c N^{-\half} \sum_{l=0}^{m-1} l^2 \lambda^l.
$$
Since $\lambda<1$ the sum is bounded independently of $m$. 
Combining the last inequality with (\ref{pre-est-inner}) completes the proof.
\end{proof}

\section{Error Analysis of the non-equispaced FFT algorithm}\label{neFFT}
The error analysis of the previous section is not complete since in
the non-equispaced FFT algorithm the complex exponential function is
approximated by the truncated Jacobi-Anger expansion. It is
therefore important to know how the highest retained order $p$ must be
selected as a function of the wave number. This will be determined in
this section.

\subsection{Error of the multivariate Jacobi-Anger approximation}
The Jacobi-Anger expansion is an expansion in Legendre polynomials.
Therefore the error of the multi-variable truncated expansion in
(\ref{trunc-ja}) is given by
\begin{equation}\label{tunc-ja2}
\hat e_k^p(t) := e_k(t) - \tilde e_k^p(t) = 
    \sum_{n=p+1}^\infty \sum_{\abs{\alpha}=n} e_\alpha L_\alpha(t)
\end{equation}
where $e_k(t)=\exp(i H \, k\cdot t)$, $\tilde e_k^p$ denotes the truncated expansion
and the coefficient $e_\alpha$ has the form
$$
e_\alpha = \int_{[-1,1]^3} e_k(t) L_\alpha(t) \,d^3t
      = { (-1)^\alpha \over \alpha! 2^\alpha} \int_{[-1,1]^3}
        (1-t^2)^\alpha \partial^\alpha e_k(t) \,d^3t.
$$
The latter form follows from the Rodrigues formula and integration by
parts. It is useful for estimating the magnitude of the coefficient
$$
\abs{e_\alpha} \leq { 1 \over \alpha! 2^\alpha} 
          \int_{[-1,1]^3} (1-t^2)^\alpha\,d^3t\, \Norm{ \partial^\alpha e_k}_\infty
            \leq  { 1 \over \alpha!\, 2^\alpha} (\pi H k)^\alpha.
$$
Since $\abs{L_n(t)} \leq 1$ for $\abs{t}\leq 1$, the truncation error
of (\ref{tunc-ja2}) can be bounded as follows
\begin{eqnarray*}
\abs{\hat e_k^p(t)} 
    &=& \sum_{n=p+1}^\infty \sum_{\abs{\alpha}=n} \abs{e_\alpha}\,, \\
    &\leq& \sum_{n=p+1}^\infty \left({\pi H\over 2} \right)^n 
              \sum_{\abs{\alpha}=n} {\abs{k^\alpha} \over \alpha!}\,,\\
    &=& \sum_{n=p+1}^\infty 
            \left({\pi H\over 2} \right)^n {\Norm{k}_1^n \over n!}\,.
\end{eqnarray*}
The last step is an application of the multivariate binomial
formula. Since the last equation is the remainder of the Taylor
expansion of the exponential function, we have the bound
\begin{equation}\label{bound-ja}
\Norm{\hat e_k^p}_{L^\infty_{[-1,1]^3}} \leq {c \over (p+1)!}  
      \left({\pi \over 4N } \Norm{k}_1 \right)^{p+1} \,. 
\end{equation}
If in (\ref{tunc-ja2}) the Taylor series instead of the Jacobi-Anger is
employed, then a very similar analysis shows that
$$
\Norm{\hat e_k^p}_{L^\infty_{[-1,1]^3}} \leq {c \over (p+1)!}  
      \left({\pi \over 2N } \Norm{k}_1 \right)^{p+1} \,.
$$
Thus the Jacobi-Anger expansion is significantly more accurate for
large values of $k$. 

\subsection{Error of the non-equispaced FFT}\label{sec-nFFT} In the
non-equispaced FFT algorithm, the kernel $G_N$ in
(\ref{fourier-series-gf}) is replaced by the kernel
\begin{equation}\label{actual-kernel}
G_N^p(x,x') = \sum_{\Norm{k}_\infty \leq N} \hat G_k
\exp\left( i\pi {k\cdot(l-l')\over N} \right) 
\tilde e_k(t_l) \tilde e_k(t'_{l'}), 
\quad x\in C_l, x'\in C_{l'}
\end{equation}
where $t_l = (x - x_l)/H$, $t'_{l'} = (x' - x_{l'})/H$. Thus the error
is given by
\begin{eqnarray*}
\lefteqn{G_N(x-x') - G_N^p(x,x') = }\\
&& \sum_{\Norm{k}_\infty \leq N} \hat G_k
\exp\left( i\pi {k\cdot(l-l')\over N} \right) 
\Big(
\hat e_k(t_l) \tilde e_k(t'_{l'}) +
\tilde e_k(t_l) \hat e_k(t'_{l'}) +
\hat e_k(t_l) \hat e_k(t'_{l'})
\Big).
\end{eqnarray*}
From (\ref{bound-ja}) and $\Norm{k}_1 \leq \sqrt{3}\Norm{k}$ the estimate
\begin{equation}\label{inter-est}
\Big|G_N(x-x') - G_N^p(x,x') \Big| \leq
{c \over (p+1)!} \left( \sqrt{3} \pi \over 4 N \right)^{p+1}
\sum_{\Norm{k}_\infty \leq N} \abs{\hat G_k} \Norm{k}^{p+1}
\end{equation}
follows. We use the identity
$$
\hat G_k = { \hat{[\Delta^2(G\chi)]}_k \over \Norm{k}^4 }
$$
and recall that from Section~\ref{sec-error} it follows that
$$
\Norm{\Delta^2(G_\delta \chi)}_{L^2} \leq c \delta^{-2}.
$$
Continuing with estimate (\ref{inter-est}) gives
\begin{eqnarray}
\lefteqn{\Big|G_N(x-x') - G_N^p(x,x') \Big|}\nonumber\\
&\leq&
{c \over (p+1)!} \left( \sqrt{3} \pi \over 4 N \right)^{p+1}
\left(\sum_{\Norm{k}_\infty \leq N} \Norm{k}^{2p-6} \right)^\half
\left(\sum_{\Norm{k}_\infty \leq N}
  \abs{\hat{[\Delta^2(G\chi)]}_k}^2\right)^\half
\nonumber\\
&\leq&
{c \over (p+1)!} \left( \sqrt{3} \pi \over 4 \right)^{p+1}
N^{-\half} \left( \delta N^2 \right)^2. \label{inter2-est}
\end{eqnarray}
We have almost completed the proof of the following theorem
\begin{theorem}
Let $\A_N^p$ the integral operator that has kernel $G_N^p$ then
$$
\big| \inner{f}{(\A_N - \A_N^p) g}\big| \leq 
{c \over (p+1)!} \left( \sqrt{3} \pi \over 4 \right)^{p+1}\!\! N^{-\half}
\,\Norm{f}_{L^2(S)} \Norm{g}_{L^2(S)}.
$$
\end{theorem}

The most important conclusion from this result is that the order $p$
in the Jacobi-Anger approximation does not have to be increased as
$N\to \infty$.

\begin{proof}
Elementary integral calculus and estimate (\ref{inter2-est}) together
with the fact that $\delta \sim N^{-2}$ imply that
\begin{eqnarray*}
\big| \inner{f}{(\A_N - \A_N^p) g}\big| &\leq& 
    \max_{x, x'\in S}\big|G_N(x-x') - G_N^p(x,x')\big|\,
    \Norm{f}_{L^1(S)} \Norm{g}_{L^1(S)}\\
 &\leq& 
  {c \over (p+1)!} \left( \sqrt{3} \pi \over 4 \right)^{p+1}\!\! N^{-\half}
    \,\Norm{f}_{L^2(S)} \Norm{g}_{L^2(S)}.
\end{eqnarray*}
\end{proof}

\section{Numerical Examples}
We have implemented the method to verify the theoretical estimates. In
this implementation the Fourier coefficients are computed
numerically. The method used for this task is completely analogous to
the computation of the Fourier coefficients of a surface density
described in Section~\ref{sec-comp-gk}. The only difference is that
the moments in (\ref{mom-surf}) are replaced by the moments of the
Green's function
$$
  m_l^\alpha(G) = 
   \int_{C_l} P_\alpha({y-x_l\over H}) G_\delta(y) \chi(y)\,d^3y
$$
where $C_l=x_l+H[-1,1]^3$. These moments are computed using
Gauss quadrature. The analysis of the error introduced by computing
$\hat G_k$ numerically parallels the discussion of
Section~\ref{sec-nFFT} and is therefore omitted. 

In the first example we compute the farfield pattern when the unit
sphere is hit with a plane wave. This is done by solving integral
equation (\ref{combined-bem}) with piecewise constant elements
combined with the spectral method. The farfield is computed from the
density using the formula
$$
\alpha(\hat x) = i\int_S \exp(-i\kappa \hat x\cdot y) (\eta + \kappa
\hat x\cdot n_y) \sigma(y) \, dS_y,\quad \hat x \in \Sphere,
$$
see, e.g.,~\cite{nedelec01}. Because of the spherical symmetry the
solution $\sigma$ as well as the farfield can be expressed in closed
form.  The coupling parameter in (\ref{combo-layer}) is
$\eta=\kappa/2$ and the linear system is solved with GMRES without any
preconditioning. Tables~\ref{tab-sphere1} and~\ref{tab-sphere1} display
the relative errors of the $L_2$-norm of $\alpha(\hat x)$ when
increasing the size of the sphere measured in wavelengths.  The
results show that the error remains bounded (actually, decreases
somewhat) when $N\sim \kappa$ and $\delta \sim 1/\kappa^2$ and are
therefore in good agreement with the theoretical estimates. The
truncation parameter of the Jacobi-Anger expansion in (\ref{tunc-ja2})
is always set to $p=4$.  In these experiments the meshwidth is
proportional to the wavelength, which is reflected in the fact that
the number of panels $n$ is quadrupled in every line. We have
implemented both filters (\ref{filter1}) and (\ref{filter2}) and set $q=5$. The
displayed results are for (\ref{filter2}), but the results for the other filter
are only marginally different.  The timings displayed are the time
per iteration and the total time, which also includes the time to
compute the Fourier coefficients $\hat G_k$. The cpu is a 3.6 ghz
Intel Xeon processor. The time per iteration
increase by a factor somewhat larger than eight when doubling the
wavenumber, which agrees well with the $\kappa^3 \log \kappa$
complexity estimate. The code stores the Fourier coefficients $\hat
G_k$ and $\hat g_k$, the moments $m^\alpha(g)$, and as well as the
orthogonal basis of the Krylov subspace generated by GMRES. For the
size of problems computed, the basis consumes the largest portion of
the overall memory usage. Since this part grows roughly like
$\kappa^2$, the growth rate of the overall storage appears slower in
the Table than the asymptotic $\kappa^3\log\kappa$ estimate.

\begin{table}[hbt]
\begin{center}
\begin{tabular}{ccccccccc}
$n$ & $N$ & size & $\delta$ & its & mem & time/itr & time & error\\
    &     & ($\lambda$) &         &     & (MB) & (sec) & (sec) &\\
\hline
5120 &   16  & 6.25 & 1.00\E{-4} & 11 &  8.0  &  0.4 &    5 & 0.076 \\
20480 &  32  & 12.5 & 2.50\E{-5} & 12 & 32.8  &  2.1 &   32 & 0.044 \\
81920 &  64  & 25   & 6.25\E{-6} & 15 & 139.8 & 18.3 &  337 & 0.039 \\
327680 & 128 & 50   & 1.56\E{-6} & 18 & 623.9 & 159.2 & 3355 & 0.034 \\
1310720 & 256 & 100 & 3.91\E{-7} & 22 & 2981.7 & 1386 & 34414 & 0.031
\end{tabular}
\caption{Results for the sphere. Lower accuracy.}
\label{tab-sphere1}
\end{center}
\end{table}

\begin{table}[hbt]
\begin{center}
\begin{tabular}{ccccccccc}
$n$ & $N$ & size & $\delta$ & its & mem & time/itr & time & error\\
    &     & ($\lambda$) &         &     & (MB) & (sec) & (sec) &\\
\hline
5120 &   16  & 3.13 & 1.00\E{-4} &  8 &  8.0  &  0.4 &    5 & 0.044 \\
20480 &  32  & 6.25 & 2.50\E{-5} & 11 & 32.8  &  2.2 &   32 & 0.021 \\
81920 &  64  & 12.5 & 6.25\E{-6} & 12 & 139.8 & 17.9 &  276 & 0.011 \\
327680 & 128 & 25  & 1.56\E{-6} & 15 & 623.9 & 159.3 & 2880 & 0.0071 \\
1310720 & 256 & 50 & 3.91\E{-7} & 18 & 2981.7 & 1378 & 28714 & 0.0051
\end{tabular}
\caption{Results for the sphere. Higher accuracy.}
\label{tab-sphere2}
\end{center}
\end{table}

\begin{figure}[hbt]
\begin{center}
\epsfig{figure=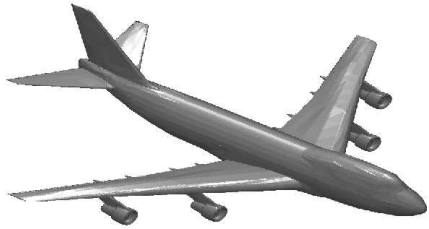,width=5in}
\end{center}
\caption{3D rendering of the airplane.}
\label{fig-boeing}
\end{figure}

To illustrate that the technique discussed in this paper can be used
for very general scatterers we include the Boeing 747 example shown in
Figure~\ref{fig-boeing}. The surface of the airplane is assumed to be
sound soft. The geometry is given by a list of vertices and triangular
panels which can be downloaded from the internet. There
are 556552 panels, and further information, such as parameterizations,
are known. We ignore the fact that there are edges and conical
vertices in the geometry and set the curvature term in
(\ref{expan-psi2}) to zero.

\begin{figure}[hbt]
\centerline{
\epsfig{figure=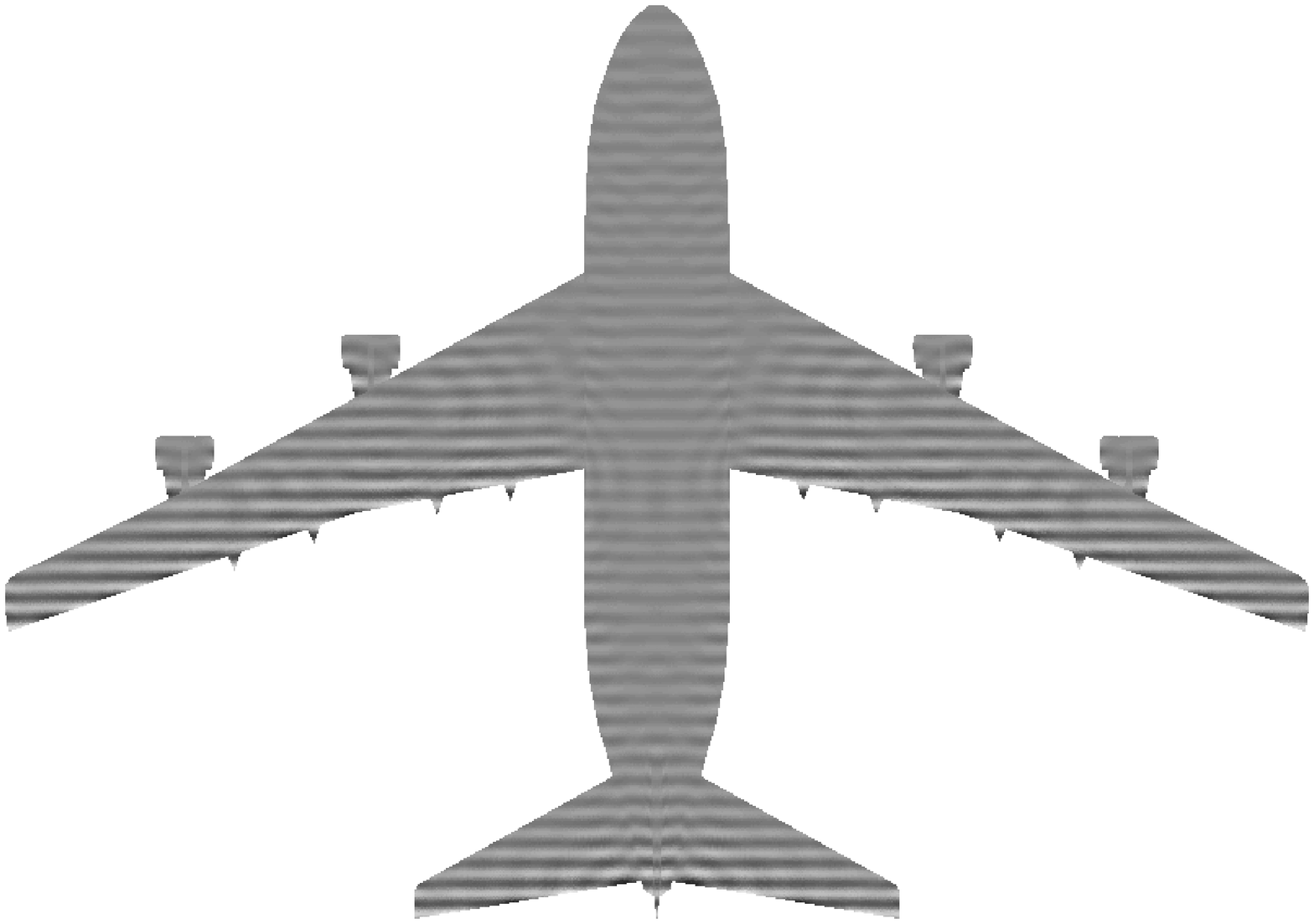,width=3in}
\epsfig{figure=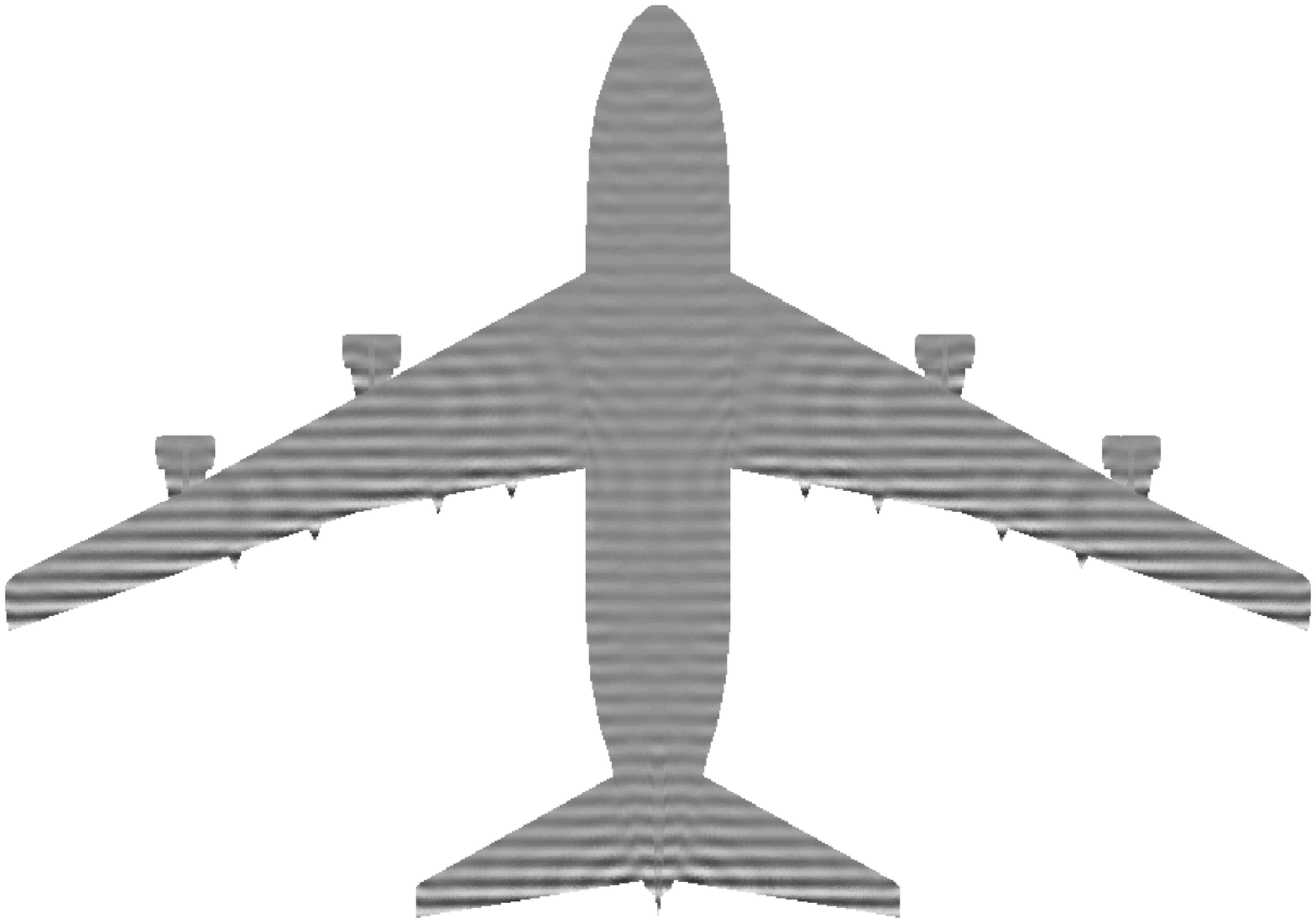,width=3in}
}
\caption{Comparison of the density (imaginary part) for $N=128$ (left) with
  $N=256$ (right); 45 wavelengths}
\label{fig-boeingLF}
\end{figure}

Figure~\ref{fig-boeingLF} compares the density for $N=128$ and $N=256$
Fourier modes. Since it is hard to spot differences in the two
solutions, it appears that already the smaller value of $N$ will give
an acceptable accuracy in many applications. The size of the scatterer
in this problem is about 45 wavelengths, the memory allocation of the
smaller problem is 906MB and the cpu time is 4355 seconds.

\begin{figure}[hbt]
\begin{center}
\epsfig{figure=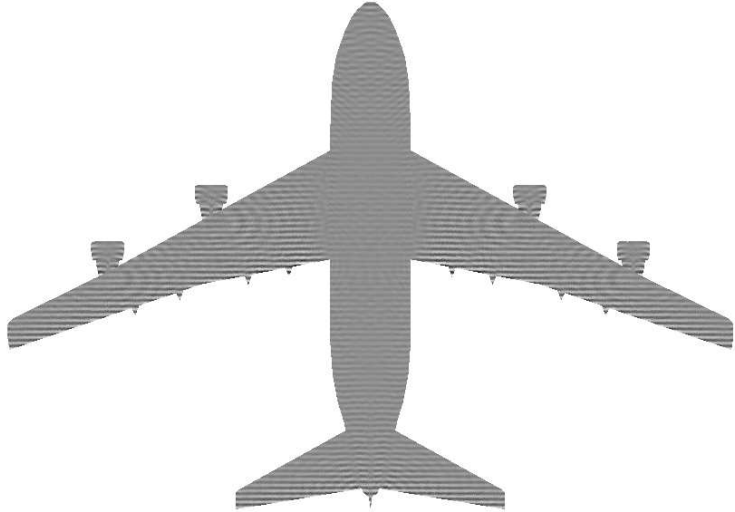,width=5in}
\end{center}
\caption{The density (imaginary part) for 90 wavelengths}
\label{fig-boeingHF}
\end{figure}

Figure~\ref{fig-boeingHF} displays the solution for 90 wavelengths and
$N=256$.  The memory allocation is 1583 MB and the cpu time is 36361
seconds. 

\section{Conclusions} We have presented a method for the computation of
scattered fields that has $O(\kappa^3 \log \kappa)$ complexity when
the meshwidth is proportional to the wavelength. Since $n\sim
\kappa^2$ the asymptotic estimate is not optimal, but because of
small constants we have been able to solve 100$\lambda$-problems in
eight to nine hours. Most of the cpu time is spent evaluating the sums
in (\ref{g-fft}) and (\ref{phi-fft}). Since this part is
embarrassingly parallel one can expect almost optimal speed up on
distributed memory multiprocessor machines. The approach generalizes
to electromagnetic scattering.

\section{Acknowledgement} The author obtained the panel description
file of the airplane from the website www.3dcafe.com. 


\end{document}